\documentclass [11pt]{amsart} 
\usepackage{epic,eepic,latexsym, amssymb, amscd, amsfonts}

\input xy
\xyoption {all}

\def \mor {{{\text{Mor}}_{d}(C, G(r,N))}}
\def \quot {{{\text{Quot}}_{d} ({\mathcal O}^N, r, C)}}
\def \gbar {{\bar g}}
\def \res {{\text {Res}}}
\def \p {{\mathbb P}}
\def \qz {{\mathsf z}}
\def \qa {{\mathsf a}}
\def \qb {{\mathsf b}}
\def \qf {{\mathsf f}}
\def \pa {{\mathsf a}}
\def \pb {{\mathsf b}}
\def \pf {{\mathsf f}}
\def \m {{\mathcal M}}
\def \n {{\mathcal N}}
\def \o {{\mathcal O}}
\def \e {{\mathcal E}}

\newtheorem{theorem}{Theorem} 
\newtheorem {lemma}{Lemma}

\newtheorem {proposition}{Proposition}

\theoremstyle{definition}
\newtheorem{remark}{Remark}

\theoremstyle {definition}

\begin{document}

\title[Intersections on the moduli space of bundles]{Counts of maps to Grassmannians and intersections on the moduli space of bundles}
\author {Alina Marian}
\address {Department of Mathematics}
\address {Yale University}
\email {alina.marian@yale.edu}
\author {Dragos Oprea}
\address {Department of Mathematics}
\address {Stanford University}
\email {oprea@math.stanford.edu}
\date{} 

\begin {abstract} We show that intersection numbers on the moduli space of stable bundles of coprime rank and degree over a
smooth complex curve can be recovered as highest-degree asymptotics in formulas of Vafa-Intriligator type. In particular, we
explicitly evaluate all intersection numbers appearing in the Verlinde formula. Our results are in agreement with previous
computations of Witten, Jeffrey-Kirwan and Liu. Moreover, we prove the vanishing of certain intersections on a suitable Quot
scheme, which can be interpreted as giving equations between counts of maps to the Grassmannian.  \end {abstract}

\maketitle

\section{Introduction}

The
Vafa-Intriligator formula counts degree $d$ maps from a smooth complex projective curve $C$ to the Grassmannian $G(r, N)$ of rank $r$ subspaces of the vector space $\mathbb C^N$, under incidence conditions with special Schubert subvarieties at fixed domain points.
When finite, the number of these maps is expressed as the sum of evaluations of a rational function at 
distinct $N^{\text{th}}$ roots of unity. Such sums are fairly complicated in general, and are not immediately expressible in a compact form. In this paper we show that intersection numbers on the moduli space $\mathcal N$ of rank $r$ and degree $d$ stable bundles on $C$, in the 
case when $r$ and $d$ are relatively prime, arise as coefficients of highest-degree terms in $N$ in Vafa-Intriligator sums. We explicitly compute these leading coefficients as iterated residues. We determine in this way all intersection numbers on $\mathcal N$ appearing in formulas of Verlinde type.  

Expressions for all top intersections of the cohomology generators on $\mathcal N$ were originally written down by Witten as
possibly divergent infinite sums over the irreducible representations of $SU(r)$ \cite {witten1}. An argument of Szenes identified the
convergent Witten sums as iterated residues \cite {Sz}. The simplest instances of these top-pairing formulas occur in rank $2$; they
received many algebro-geometric or symplectic proofs. By contrast, the higher rank case is more involved. A symplectic derivation of all
of the intersection theory of $\mathcal N$ was given by Jeffrey and Kirwan \cite {jeffreykirwan}, with intersection numbers expressed as
iterated residues. The Jeffrey-Kirwan approach is very powerful, yet not easily accessible to the uninitiated: technical difficulties
concerning localization on singular non-compact moduli spaces have to be carefully considered. A different approach is due to Liu who
obtained, by heat kernel methods, some of the intersection numbers as Witten sums \cite {liu}. Another derivation via large level limits
of Verlinde-type formulas was hinted at in the last section of \cite {TW}. Finally, results for arbitrary groups were obtained by Meinrenken \cite{M} via symplectic geometry.

The goal of our paper is to present a surprisingly simple and entirely finite dimensional algebro-geometric derivation of the Witten-Szenes-Jeffrey-Kirwan-residue formulas. The inspiration for the current work was the connection between the intersection theory of $\mathcal N$ and the enumerative geometry of Grassmannians, which emerged several years ago in the mathematics and physics literature \cite {bertram1} \cite{witten2}. On the physics side, Witten equates the Verlinde formula giving the dimension of non-abelian (possibly parabolic) theta functions on $C$ with a certain count of maps to a Grassmannian. A precise statement (in the absence of parabolic
structures) is given later in this introduction.  Witten argues that the Grassmannian nonlinear sigma model and the gauged WZW model of $U(r)/U(r)$, obtained after integrating out matter fields, arise as different large-distance limits of the same physical theory, a gauged linear  sigma model. Some of the topological amplitudes of the gauged WZW model, which has a direct link to the Verlinde algebra, and of the Grassmannian sigma model should therefore be closely related. In genus $0$, this relation was mathematically established by Agnihotri \cite {A}, who proved that the small quantum cohomology of the Grassmannian $G(r,N)$ is isomorphic to the Verlinde algebra of the unitary group $U(r)$ at level $N-r$. In higher genus, we will extract (roughly half of) the intersection theory of the moduli space of bundles from the count of maps to the Grassmannian. In the limit $N\to \infty$, it may be possible to relate our method (even if only philosophically) to the more abstract approach pursued in \cite {T}\cite{TW}. Note however that the bounded geometry of our setup lends itself to direct intersection-theoretic calculations, via simpler machinery, avoiding the use of {\it $K$-theory} on the moduli {\it stack}.\\

To start, we let $\mathcal N$ be the moduli space of rank $r$, degree $d$ stable bundles with fixed determinant on a smooth complex projective curve $C$ of genus $g\geq 2$. We will frequently use the notation $\gbar=g-1.$ We assume throughout that $r$ and $d$ are coprime. Let $\mathcal V$ be the universal bundle on $\mathcal N\times C$. We K\"{u}nneth-decompose its Chern classes $$c_i(\mathcal V)= a_i\otimes 1
+ \sum_{j=1}^{2g} b_i^{j}\otimes \delta_j+f_i\otimes \omega, \,\,\, 1\leq i\leq r,$$ with respect to a fixed symplectic basis $1,
\delta_1, \ldots, \delta_{2g}, \omega$ of $H^{\star}(C)$. The classes $a_i$, $b_i^{j}$, $f_i$ thus obtained are generators for the
cohomology of $\mathcal N$ \cite{atiyahbott}. They are however not canonical, since the universal bundle $\mathcal V$ is only defined up to tensoring with a line bundle from the base $\mathcal N$.

A normalization process is required to kill this ambiguity and obtain {\it invariant} classes. To this end, we write $x_1, \ldots, x_r$
for the Chern roots of $\mathcal V$, and let \begin{equation}\label{normalized}\bar x=\frac{1}{r}\left(x_1+\ldots+x_r\right), \,\,\, \bar
x_i = x_i - \bar x, \,1\leq i\leq r.\end{equation} The $i^{\text {th}}$ symmetric elementary polynomial in the normalized variables is
written as $$s_i(\bar x_1, \ldots, \bar x_r)=\bar a_i \otimes 1+ \sum_{j=1}^{2g} \bar b_i^{j}\otimes \delta_j + \bar f_i \otimes \omega,
\,\,\, 1\leq i\leq r.$$ For instance, $$\bar f_2=f_2-\frac{d(r-1)}{r}a_1=\frac{1}{2r}c_1(\mathcal N),$$ $$\bar a_1=0, \text{ and }\,\,
\bar a_i=s_i(\bar \theta_1, \ldots, \bar \theta_r).$$ Here $\theta_1, \ldots, \theta_r$ are the Chern roots of $\mathcal V|_{\mathcal N \times
\text{point}}.$ \\

We can now explain the main result. Consider the formal variables $\mathsf X_1,\ldots, \mathsf X_r$ and $\mathsf Y_1, \ldots, \mathsf
Y_{r-1}$ such that \begin{equation}\label{xy}\mathsf X_1+\ldots+\mathsf X_r=0,\end{equation} \begin{equation}\label{xy1}\mathsf
Y_1=\mathsf X_1-\mathsf X_2, \ldots, \mathsf Y_{r-1}=\mathsf X_{r-1}-\mathsf X_r.\end{equation} 
Define moreover the function
\begin{equation}\label{m}\mathsf L=\left\{\frac{d}{r}\right\}\mathsf
Y_1+\ldots+\left\{\frac{d(r-1)}{r}\right\}\mathsf Y_{r-1},\end{equation} where $\left\{\right\}$ denotes the fractional part.

\begin{theorem}\label{main} For any polynomial $\mathsf P$ in the normalized $\bar a$ classes, we have \begin{eqnarray}\label{mainformula}\int_{\mathcal
N}\exp({\bar f_2})\, \mathsf P(\bar a_2, \ldots, \bar a_r) &=&(-1)^{\gbar\binom{r}{2}}\,r^{\gbar} \,\res_{\mathsf
Y_{1}=0}\ldots \res_{\mathsf Y_{r-1}=0} \frac{1}{e^{\mathsf Y_1}-1}\cdots\frac{1}{e^{\mathsf Y_{r-1}}-1}\times \nonumber \\ &\times& \exp\left(\mathsf L\right)\cdot\frac{\mathsf P(s_2(\mathsf X), \ldots, \mathsf s_r(\mathsf X))}{\prod_{i<j} (\mathsf X_i-\mathsf X_j)^{2\gbar}}.\end{eqnarray} Here $s_k$ are the elementary symmetric functions. The iterated residue is computed from right to left,
at each step keeping all but one of the $\mathsf Y$s fixed.  \end{theorem}

Theorem \ref{main} will be obtained by exploiting the connection between the intersection theory of $\mathcal N$ and that of a
suitable compactification of the scheme $\mor$ of degree $d$ morphisms from $C$ to the Grassmannian $G(r,N).$ This 
compactification, the scheme $\quot$ constructed by Grothendieck,
is the fine moduli space of short exact sequences \begin{equation}\label{seq}0\to E\to \mathcal O^{N}\to F\to 0,\end{equation}
where $F$ is a degree $d$, rank $r$ quotient {\it sheaf} of the trivial bundle. With the aid of the dual universal bundle
$\mathcal E^{\vee}$ on $\quot \times C$, we obtain, via the K\"{u}nneth decomposition, cohomology classes denoted by $\qa, \qb,
\qf$: \begin{equation}\label{ekunneth}c_i(\mathcal E^{\vee})=\qa_i\otimes 1+\sum_{j=1}^{2g} \qb_i^{j}\otimes \delta_j +
\qf_i,\,\,\, 1\leq i\leq r.\end{equation} When the degree $d$ is large, it was established in \cite{bertram1} that $\quot$ is
irreducible, generically smooth of the expected dimension $$e=\chi(E^{\vee}\otimes F)=Nd-r(N-r)\gbar.$$ In general, this space
may be singular and of the wrong dimension. Nonetheless, $\quot$ comes equipped with a canonical perfect obstruction theory, and
hence with a virtual fundamental cycle \cite{vi} $$\left[\quot\right]^{vir}\in A_{e}(\quot).$$ Polynomials in the $\qa, \qb,
\qf$ classes can then be evaluated against this virtual fundamental cycle. In particular, the Vafa-Intriligator formula
\eqref{vif} proved in \cite{sieberttian}, \cite{vi} expresses top intersections of $\qa$-classes as sums over $N^{\text{th}}$
roots of unity. \\

Note that all intersection numbers appearing in the Verlinde formula are covered by Theorem \ref{main}. A derivation of this formula {\it via Riemann-Roch} was obtained in the last section of \cite {jeffreykirwan}; the computation for arbitrary structure groups was pursued in \cite{bismut}. Denoting by $\mathcal
L$ the ample generator of $\text{Pic }(\mathcal N)$, we have $c_1(\mathcal L)=r\bar {f_2}$, hence \begin{equation}
\label{verlinde0}\chi(\mathcal L^{s})=\int_{\mathcal N} \exp(sr\bar f_2) \,\text{Todd}(\mathcal N)=\int_{\mathcal N}\exp\left((s+1)r\bar
f_2\right) \hat A(\mathcal N).\end{equation} Here, $\hat A(\mathcal N)$ is a polynomial in the $\bar a_i$ classes, which is explicitly
given in terms of the Chern roots $\theta_i$ as $$\hat A(\mathcal
N)=\prod_{i<j}\left(\frac{\theta_i-\theta_j}{2\sinh\frac{\theta_i-\theta_j}{2}}\right)^{2\gbar}.$$ The evaluation of the Verlinde Euler
characteristic \eqref{verlinde0} is immediate from equation \eqref{mainformula}. The emerging residue answer can be recast effortlessly as an intersection number on $\text{Quot}_d (\mathcal O^{r(s+1)}, r, C)$ by a {\it backwards} application of the residue formula \eqref{intresa} of section \ref{residues}. The result is the following striking equality, derived by physical considerations in \cite {witten2},
\begin{equation}\label{verlinde}\chi(\mathcal L^s)=\frac{1}{(s+1)^{g}} \int_{\left[\text {Quot}_d (\mathcal O^{r(s+1)}, r, C)\right]^{vir}}
\qa_r^{s(d-r\gbar)+d}.\end{equation} Moreover, when $d$ is large, intersections of $\qa$ classes on the Quot scheme have enumerative
meaning \cite {bertram}. In particular, the right-hand-side integral of \eqref{verlinde} is the finite count of degree $d$ maps to the
Grassmannian $G(r, r(s+1))$ with incidences at fixed $s(d-r\gbar)+d$ domain points with sub-Grassmannians $G(r, r(s+1)-1)\hookrightarrow
G(r,r(s+1))$ in general position.\\

In a different direction, the argument which gives Theorem \ref{main}, combined with a rescaling trick we learned
from \cite {EK}, leads to a vanishing result about intersections on Quot, which will be presented in the last section of this work. \\
  
The paper is organized as follows. First, we explain the setup of \cite {marian} which relates intersections on the moduli space of stable bundles to intersections on $\quot$ in the large $N$ regime. Using as starting point the Vafa-Intriligator formula, the relevant evaluations on $\quot$
are cast as iterated residues in Section \ref{residues}. Upon extracting the appropriate asymptotic coefficients of these residues, we immediately obtain in Section \ref{der} the formulas of Theorem \ref{main}. For completeness, we also indicate how the results are expressed as infinite sums indexed by the irreducible representations of $SU(r)$, as in \cite{witten1}. Finally, the vanishing of intersections on $\quot$ will be proved in the last section.

\section {Verlinde-type intersections through the Quot scheme}

We start by explaining the setup for the proof of Theorem \ref{main}. Traditionally, one studies the intersection theory on 
the moduli space $\mathcal N$ of rank $r$ degree $d$ stable bundles with fixed determinant. However, our computations will be most naturally carried on the moduli space $\mathcal M$ of stable bundles $V$ with varying 
determinant. These intersections will be transfered to $\mathcal N$ via the the degree $r^{2g}$ {\'e}tale cover $$\tau: \mathcal N\times J\to \mathcal M,$$ 
given by tensoring with degree $0$ line bundles in the Jacobian $J$ of $C$. When switching from $\mathcal N$ to $\mathcal M$, we will abusively use the 
same 
notations for the universal bundle, its Chern roots and the K\"{u}nneth components of its Chern classes.

The main technical ingredient of our argument is a precise relationship between the intersection theory of $\mathcal M$ and that
of $\quot$. We will assume that $d$ is large compared to $N$, $r$ and $g$ to ensure that $\quot$ is irreducible of the expected
dimension. We will think of the points of $\quot$, that is of the short exact sequences \eqref{seq}, as $N$ tuples of sections of
the bundle $E^{\vee}$ which generically generate the fiber. Requiring that $V=E^{\vee}$ be stable and not demanding that the
sections generically generate, we arrive at a different moduli space, birational to $\quot$, which we denote by $\p_{N,r,d}$.
$\p_{N, r, d}$ finely parametrizes pairs $(V, \phi)$ where $V$ is a stable rank $r$ degree $d$ bundle on $C$, and $\phi$ is a
non-zero $N$ tuple of holomorphic sections, considered projectively $$\phi: \mathcal O^N\to V.$$ When $d$ is large, the space
$\pi:\p_{N, r, d}\to \mathcal M$ is the projective bundle $\p_{N, r, d}=\p(\mathcal H^{N}),$ where $$\mathcal H=pr_{\star}
\mathcal V,$$ with $pr$ denoting the projection from $\m \times C$ to $\m$.  

There is a universal morphism $$\Phi: \o^N \rightarrow {\mathcal U} \, \, \, \text{on} \, \, \, \p_{N, r, d} \times C,$$ 
and it is easy to see that \begin{equation}\label{uversv}\mathcal U=\pi^{\star} \mathcal 
V\otimes \mathcal O_{\p}(1).\end{equation} We K\"{u}nneth decompose the Chern classes of $\mathcal U$ as \begin{equation}\label{ukunneth} c_i(\mathcal 
U)=\pa_i\otimes 1+ \sum_{j=1}^{2g} \pb_i^{j} \otimes \delta_j + \pf_i\otimes \omega, \,\, 1\leq i\leq r.\end{equation} 

The reader who compared \eqref{ekunneth} and \eqref{ukunneth} may have noticed the abusive notation. The reason for using the same symbols for analogous, 
but
certainly different classes on $\p_{N, r, d}$ and $\quot$, is that both moduli spaces and their universal structures agree on an open subscheme. Usually we will
carefully distinguish between the different $\qa, \qb, \qf$ classes, by always mentioning the moduli spaces where the intersections are computed. As
before, we let $\bar \pa_i, \bar \pb_i^{j}, \bar \pf_i$ be the normalized classes either on $\quot$ or $\p_{N, r, d}$. (Note however that the universal 
structures on $\quot$ and $\p_{N, r, d}$ {\em{are}} canonical.)

We summarize our setup in the following diagram \begin{center}$\xymatrix{\mor \hookrightarrow \quot \ar@{<.>}[r] & \p_{N, r, d}\ar[d]^{\pi}& \\ & \mathcal M & \mathcal N\times J\ar[l]^{\tau}}.$ \end{center}

Let us consider the intersection product 
\begin{equation}
\label{mainquot}
\int_{\quot} {\mathsf P} (\bar \qa_2, \ldots, \bar \qa_r) \cdot \qa_r^M,
\end{equation}
where
$\mathsf P (\bar \qa_2, \ldots, \bar \qa_r)$ is a polynomial in the normalized $\bar {\mathsf a}$ 
classes of
total weighted degree at most $r^2\gbar+1$, and $M$ is such that \begin{equation}\label{M}\text {deg }\mathsf P+ r M=
Nd-r(N-r)\gbar.\end{equation} This choice is possible since $d$ and $r$ are coprime. 
In addition to our initial assumption that $d$
is large compared to $N$, $r$ and $g$, we will require that $N$ be large with respect to $r$ and $g$. This warrants the transfer
of certain intersections from $\quot$ to $\p_{N, r, d},$ by making sure that the nonoverlapping loci of $\quot$ and $\p_{N,
r, d}$ are avoided. Indeed,
when $N$ is large, the main theorem of \cite 
{marian} asserts that the class ${\mathsf P} (\bar \qa_2, \ldots, \bar \qa_r) \cdot \qa_r^M$ evaluates identically on $\quot$ and
${\p_{N, r, d}}.$ It is moreover easy to observe, using \eqref{uversv}, that the
normalized $\bar a$ classes on $\p_{N, r, d}$ and $\mathcal M$ are related by pullback \begin{equation}\label{pullback}\mathsf
P(\bar \pa_2, \ldots, \bar \pa_r)=\pi^{\star} \mathsf P(\bar a_2, \ldots, \bar a_r).\end{equation}
Consequently, we have that for $N$ large relative to $r$ and $g$, and $d$ large relative to 
$N$, $r$ and $g$, 
\begin{equation}\label{alina} \int_{\quot} \mathsf P (\bar \qa_2, \ldots, \bar \qa_r) \cdot \qa_r^{M} = \int_{\p_{N, r, d}}
\mathsf P (\bar \pa_2, \ldots, \bar \pa_r) \cdot \pa_r^{M} = \int_{\mathcal M} {\mathsf P(\bar a_2, \ldots, \bar 
a_r)} \cdot \pi_{\star} (\pa_r^{M}).\end{equation} 

The observation which lies at the heart of the argument for Theorem \ref{main} is that the top-degree term in $N$ of the 
intersection \eqref{mainquot}
is closely related to the intersection numbers which are the subject of the theorem.  
We will therefore study the leading behavior in $N$ of \eqref{alina}. $N$ is however assumed to be small compared to the
degree $d$. We remove this assumption by observing that all intersections appearing in \eqref{alina} depend only on $d$ modulo $r$. This is
clear for the rightmost intersection by the arguments below, and it follows for the one on $\quot$ from the explicit evaluation via the
Vafa-Intriligator formula. Alternatively, we can reduce the degree $d$ on $\quot$ directly, making use of the virtual
fundamental class on $\quot$ constructed in \cite {vi} and Theorem $2$ there, which compares the virtual cycles for same
values of $d$ modulo $r$.\\

Let us examine the pushforward $\pi_{\star} (\pa_r^{M}).$ By \eqref{uversv}, the Chern roots of $\mathcal U$ restricted
to $\p_{N, r, d}\times \{\text{point}\}$ equal $\pi^{\star}\theta_1+\zeta, \ldots, \pi^{\star} \theta_r+\zeta$, with $\theta_1, \ldots,
\theta_r$ being the Chern roots of $\mathcal V|_{\mathcal M\times {\text{point}}}$, and $\zeta=c_1(\mathcal O_{\p}(1)).$ Hence, \begin{eqnarray}\label{segre}\pi_{\star}
(a_r^{M})&=&\pi_{\star}\left((\pi^{\star}\theta_1+\zeta)^M\cdot\ldots\cdot (\pi^{\star}\theta_r+\zeta)^{M}\right) \\ &=&\sum_{l_1,
\ldots, l_r} \binom{M}{l_1}\ldots \binom{M}{l_r}\,\theta_1^{l_1}\cdot \ldots\cdot \theta_r^{l_r} \cdot\pi_{\star}
(\zeta^{rM-|l|})\nonumber\\&=&\sum_{l_1, \ldots, l_r} \binom{M}{l_1}\ldots \binom{M}{l_r}\,\theta_1^{l_1}\cdot \ldots\cdot \theta_r^{l_r}
\cdot s_{k}\left(\mathcal H^{N}\right)\nonumber,\end{eqnarray} where $$k=rM-|l|-\left(N \text{rank } \mathcal
H-1\right)=r^2\gbar+1-|l|-\text {deg } \mathsf P.$$ The following lemma helps determine the leading term in $N$ of the above expression.

\begin{lemma}\label{segreasymp} As $N\to \infty$, we have $$s_k(\mathcal H^N)=\frac{(-1)^k}{k!} c_1(\mathcal
H)^{k}\cdot N^{k} +\text { lower order terms in } N.$$ \end{lemma} \noindent{\bf Proof.} Letting $h_1\ldots h_s$ be the Chern roots of $\mathcal H$, we
have \begin{eqnarray} s_{k}(\mathcal H^N)&=& \left[\frac{1}{c(\mathcal H)^{N}}\right]_{(k)}=\left[\frac{1}{(1+h_1)^N\ldots
(1+h_s)^N}\right]_{(k)}=\nonumber \\ &=&\sum_{k_1+ \ldots+ k_s=k}\binom{-N}{k_1}h_1^{k_1}\cdot \ldots \cdot
\binom{-N}{k_s}h_s^{k_s} =\nonumber\\ &=&\sum_{k_1+\ldots+ k_s=k}\frac{(- N)^{k_1}}{k_1!}h_1^{k_1}\cdot \ldots \cdot
\frac{(-N)^{k_s}}{k_s!}h_s^{k_s}+\text{lower order terms in } N=\nonumber\\ &=& (-1)^k  \frac{(h_1+\ldots+h_s)^{k}}{k!} \cdot N^k+\text{ 
lower order terms in } N= \nonumber \\ & = & \frac{(-1)^k}{k!} c_1(\mathcal H)^k \cdot N^k+\text {lower order terms in } 
N.\nonumber\end{eqnarray}\qed

Observing that $$\binom{M}{l}=\frac{1}{l!} \left(\frac{d}{r}-\gbar\right)^{l}N^{l}+\text {lower order terms in } N,$$ and using Lemma
$\ref{segreasymp}$, we compute that the highest exponent of $N$ in \eqref{segre} is $$\mathtt e({\mathsf P})=r^2\gbar+1-\text {deg
}\mathsf P,$$ and the corresponding coefficient equals \begin{eqnarray}\label{coeff}&\sum_{l_1, \ldots,
l_r}&\frac{\theta_1^{l_1}}{l_1!}\left(\frac{d}{r}-\gbar\right)^{l_1}\cdot \ldots\cdot
\frac{\theta_r^{l_r}}{l_r!}\left(\frac{d}{r}-\gbar\right)^{l_r}\frac{(-c_1(\mathcal H))^{k}}{k!}=\\&=&
\left[\exp\left(\left(\frac{d}{r}-\gbar\right)(\theta_1+\ldots+\theta_r)-c_1(\mathcal H)\right)\right]_{(\mathtt e({\mathsf
P}))}\nonumber\\&=&\left[\exp\left(\left(\frac{d}{r}-\gbar\right)a_1-\left((d-\bar g)a_1-\sum_{j=1}^{g}b_1^{j}
b_1^{j+g}-f_2\right)\right)\right]_{(\mathtt e({\mathsf P}))}\nonumber\\&=&\left[\exp\left (f_2 -\frac{d(r-1)}{r} a_1 + \sum_{j=1}^{g} 
b_1^{j}
b_1^{j+g}\right)\right]_{(\mathtt e({\mathsf P}))}.\nonumber\end{eqnarray} In the third line, $c_1(\mathcal H)$ is computed by
Grothendieck-Riemann-Roch for the projection $pr: \m \times C \rightarrow \m.$ As a consequence of \eqref{pullback}, \eqref{segre}, and
\eqref{coeff}, we obtain that \eqref{alina} equals $$N^{\mathtt e({\mathsf P})} \int_{\mathcal M} \exp\left(
f_2 - \frac{d(r-1)}{r} a_1 + \sum_{j=1}^{g} b_1^{j} b_1^{j+g}\right)\cdot \mathsf P(\bar a_2, \ldots, \bar a_r) +\text{lower order terms 
in}\, \, N.$$ 

We pull back this expression under $\tau: \n \times J \rightarrow \m$ observing that, since $\mathcal N$ is simply connected
\cite{atiyahbott}, we have $$\tau^{\star} \left(\sum_{j=1}^{g} b_1^{j} b_1^{j+g}\right)=r^2\Theta, \text{ and } \tau^{\star} \left (f_2 
- \frac{d(r-1)}{r} a_1 \right ) =
\bar f_2 - r(r-1)\Theta.$$ Here $\Theta$ is the theta class on the Jacobian $J$. We find that the highest coefficient above is
\begin{equation}\nonumber\label{coeff2}\frac{1}{r^{2g}}\int_{\mathcal N\times J} \exp\left(\bar f_2+r\Theta\right) \mathsf P(\bar a_2,
\ldots, \bar a_r) =\frac{1}{r^{2g}}\int_{\mathcal N} \exp(\bar f_2)\,\mathsf P(\bar a_2, \ldots, \bar a_r)\cdot
\int_{J}\exp(r\Theta)=\end{equation} $$=\frac{1}{r^g}\int_{\mathcal N} \exp(\bar f_2)\, \mathsf P(\bar a_2, \ldots, \bar a_r).$$ We
therefore conclude that \begin{equation}\frac{1}{r^g}\label{finalint}\int_{\mathcal N} \exp(\bar f_2)\, \mathsf P(\bar a_2, \ldots, \bar
a_r) =\left[N^{r^2 \gbar+1-\text {deg } \mathsf P}\right] \int_{\quot} \mathsf P (\bar \qa_2, \ldots, \bar \qa_r) \cdot \qa_r^{M},
\end{equation} where the brackets denote taking the coefficient of the given power of $N$.

\section{Intersections on the Quot scheme as residues} \label{residues}
In order to calculate the right-hand side of \eqref{finalint} and thus complete the proof of Theorem \ref{main}, we will make use of the Vafa-Intriligator formula 
\cite{sieberttian}, \cite{vi}, and will rewrite the resulting expression as a residue. It will be then easy to extract the leading $N$ 
coefficient in the next section. 
Recasting Vafa-Intriligator as a residue is essentially done in \cite{Sz} and \cite{jeffreykirwan}, but for completeness we 
would like to indicate the argument.

To start, we recall the statement of the Vafa-Intriligator formula.\\ 

{\it Let $\mathsf A$ be a symmetric top-degree polynomial in the Chern roots $\qz_1, \ldots, \qz_r$ of the dual universal sheaf $\e^{\vee}$ restricted to $\quot \times \{\text {point}\}.$ Then
\begin{equation}\label{vif}
\int_{[\quot]^{vir}} \mathsf A (\qz_1, \ldots, \qz_r) = \mathsf u N^{r\gbar} \sum_{\lambda_1, \ldots, \lambda_r} \frac{\mathsf A (\lambda_1, \ldots 
\lambda_r 
)\cdot(\lambda_1 \cdots \lambda_r)^{-\gbar}}{\prod_{i<j} (\lambda_i - \lambda_j )^{2\gbar}}.
\end{equation}
Here $(\lambda_1, \ldots, \lambda_r)$ are {\it ordered} tuples of distinct $N^{\text{th}}$ roots of $1$ and $\mathsf u=(-1)^{\gbar\binom{r}{2}+d(r-1)}.$ }\\

The intersection number on the right-hand side of \eqref{finalint} is therefore
\begin{equation} \label{vinow}
\int_{\quot} \mathsf P (\bar \qa_2, \ldots, \bar \qa_r) \cdot \qa_r^{M}= \mathsf u N^{r\gbar}\cdot
\sum_{\lambda_1,\ldots, \lambda_r} \mathsf Q(\lambda_1, \ldots, \lambda_r) \cdot
\frac{(\lambda_1\ldots\lambda_r)^{M-\gbar}}{\prod_{i<j}\left(\lambda_i-\lambda_j\right)^{2\gbar}}
.\end{equation}
The
polynomial $\mathsf Q$ is obtained from $\mathsf P$ using the symmetric functions in the normalized variables $\bar x_1, \ldots, \bar x_r$ defined as in
\eqref{normalized}, \begin{equation}\label{Q}\mathsf Q(x_1, \ldots, x_r)=\mathsf P(s_2(\bar x_1, \ldots, \bar x_r), \ldots, s_r(\bar x_1, \ldots, \bar
x_r)).\end{equation}

Note now that by \eqref{M}, $$\deg \mathsf Q - r(r-1) \gbar + r (M - \gbar) \equiv 0 \mod N, $$ so each summand of \eqref{vinow} is invariant under rescaling 
by
an $N^{\text{th}}$ root of unity. After rescaling the roots of unity $\lambda_1, \ldots, \lambda_r$ by $\lambda_r^{-1}$, we may assume $\lambda_r=1$. This normalization changes the power of $N$ by $1$, to account for each possible value of $\lambda_r$. Furthermore, we will allow {\it unordered} tuples in the sum
\eqref{vinow} at the expense of the prefactor $\frac{1}{r!}$. Setting \begin{equation}\label{R}\mathsf R(x_1, \ldots, x_r) = \mathsf Q(x_1, \ldots, x_r)
\cdot \prod_{i<j} (x_i-x_j)^{-2\gbar}.\end{equation} we rewrite the right hand side of \eqref{vinow} as
\begin{equation}\label{suminr} N^{r\gbar +1}  \frac{\mathsf u}{r!}\sum_{\lambda_1, \ldots, \lambda_{r-1}} \mathsf
R(\lambda_1, \ldots, \lambda_{r-1}, 1) \cdot (\lambda_1 \cdots \lambda_{r-1})^{M-\gbar},\end{equation} where we sum over tuples of roots of $1$ which are
pairwise distinct and not equal to $1$.

Further set \begin{equation}\label{alpha}\alpha_1 = \frac{\lambda_1}{\lambda_{2}},\,\,\, \ldots\,
\alpha_{r-2}=\frac{\lambda_{r-2}}{\lambda_{r-1}},\,\,\, \alpha_{r-1}=\lambda_{r-1},\end{equation} so that $$\lambda_1=\alpha_1\cdot 
\alpha_2 \cdot\ldots\cdot
\alpha_{r-1},\,\,\ldots\, ,\lambda_i=\alpha_i\cdot\ldots\cdot \alpha_{r-1},\,\,\,\ldots\, ,\lambda_{r-1}=\alpha_{r-1}.$$ The sum in \eqref{suminr}
becomes \begin{equation} \label{sumalpha} \sum_{(\alpha_1, \ldots,\alpha_{r-1}) \in {\mathcal C}} {\mathsf R} (\alpha_1 \cdot
\ldots\cdot\alpha_{r-1} , \ldots, \alpha_{r-1} ,1) \cdot \left(\alpha_1 \cdot{\alpha}_2^{2}\cdot\cdots
\cdot\alpha_{r-1}^{r-1}\right)^{M-\gbar}.\end{equation} Here $\mathcal C$ is the set of all (arbitrary length) tuples $(\alpha_1, \ldots,
\alpha_{s})$ of $N^{\text{th}}$ roots of $1$ such that all products $$\alpha_i\cdot \alpha_{i+1} \cdots \alpha_j\neq 1,\,\,\, 
\text {for }
1\leq i\leq j\leq s.$$ Moreover, let us define the integers $0\leq m_i<N$ such that \begin{equation}\label{reducedexp}i(M-\gbar)\equiv
m_i \mod N, \, \, \, \, 1\leq i \leq r-1.\end{equation} Then, we can rewrite \eqref{sumalpha} as 
\begin{equation}\label{reducedexpsum}\sum_{(\alpha_1, \ldots,
\alpha_{r-1}) \in {\mathcal C}} {\mathsf R} (\alpha_1 \cdot \ldots\cdot\alpha_{r-1} , \ldots, \alpha_{r-1} ,1) \cdot
\alpha_1^{m_1}\cdots\alpha_{r-1}^{m_{r-1}}.\end{equation} Observe that we have \begin{equation}\label{min}\frac{m_i}{N}\to \left\{\frac{di}{r}\right\} \text { as } N\to \infty.\end{equation}

Consider now the meromorphic one-form in the variables $y_1, \ldots, y_{r-1}$
$$\Omega_1 = \frac{dy_{r-1}}{y_{r-1}} \cdot \frac{N}{y_{r-1}^N -1}\cdot{\mathsf R} (y_1 \cdot\ldots\cdot y_{r-1} , \ldots, 
y_{r-1} , 1) \cdot y_1^{m_1}\cdots y_{r-1}^{m_{r-1}}.$$ We think of the variables $y_1, \ldots , y_{r-1}$ as standing in a relation to the original $x_1, \ldots, x_r$ used in \eqref{R} completely similar to 
the relation that the $\alpha$s bear to the original $\lambda$s. That is,  
\begin{equation} \label{xs} \frac{x_i}{x_r} = y_i \cdot \cdots \cdot y_{r-1}, \, \, \, 1\leq i \leq r-1. \end{equation}

One checks that for each fixed $(\alpha_1,\ldots, \alpha_{r-2})\in \mathcal C$, the form $\Omega_1 (\alpha_1, \ldots, \alpha_{r-2}, y_{r-1})$ is 
meromorphic on the projective line, with poles at all the $N^{\text{th}}$ roots of unity, and nowhere else. This step makes essential use of the fact
that the exponents $m_i$ of the $y_i$ are roughly non-zero subunitary fractions of $N$, thus eliminating the possibility of 
poles at $0$ or $\infty.$ By the global residue
theorem the sum of the residues of $\Omega_1$ at these poles is zero.

Note that we have simple poles at the roots $\nu$ such that $(\alpha_1, \ldots, \alpha_{r-2}, \nu ) \in {\mathcal C}$ {i.e.,} such
that $$\nu\neq 1,\,\,\nu\neq(\alpha_i\cdots\alpha_{r-2})^{-1},\,\, 1\leq i\leq r-2.$$ Their residues are $$\res_{y_{r-1}=\nu}\, \Omega_1
(\alpha_1, \ldots, \alpha_{r-2}, y_{r-1})={\mathsf R} (\alpha_1 \cdots \alpha_{r-2}\, \nu , \ldots, \alpha_{r-2}\, \nu ,\, \nu, 1) \cdot
\alpha_1^{m_1}\cdots \alpha_{r-2}^{m_{r-2}}\cdot \nu^{m_{r-1}}.$$ As a consequence, the sum in \eqref{reducedexpsum} is
\begin{eqnarray}&&\sum_{(\alpha_1, \ldots, \alpha_{r-1}) \in {\mathcal C}} {\mathsf R} (\alpha_1 \cdot \ldots\cdot\alpha_{r-1} , \ldots,
\alpha_{r-1} ,1) \cdot \alpha_1^{m_1}\cdots\alpha_{r-1}^{m_{r-1}}=\nonumber\\ &=& - \sum_{(\alpha_1, \ldots, \alpha_{r-2}) \in \mathcal
C} \left({\text{Res}}_{y_{r-1} = 1}\, \Omega_1 +\sum_{i=1}^{r-2} {\text{Res}}_{y_{r-1} = (\alpha_i \cdots \alpha_{r-2})^{-1}}
\,\Omega_1\right) {\Bigg|}_{ {\stackrel{y_1=\alpha_1}{y_{r-2}=\alpha_{r-2}}}}.\nonumber\end{eqnarray}

For each $i < r-2$ in the above sum of residues, define the
rescaled variables $$\widetilde y_{r-1}=\alpha_i \cdot \ldots \cdot \alpha_{r-2}\cdot y_{r-1},$$ $$\widetilde \alpha_{r-2}=\left(\alpha_i \cdot \ldots \cdot \alpha_{r-2}\right)^{-1} \cdot\alpha_{r-2},$$ $$\widetilde \alpha_i=\left(\alpha_{i} \cdot \ldots \cdot \alpha_{r-2}\right)^{-1} \cdot\alpha_i$$ $$\widetilde \alpha_{i-1}
=\left (\alpha_{i} \cdot \ldots
\cdot \alpha_{r-2} \right )\cdot \alpha_{i-1}.$$ We keep the rest of the $\alpha$s unchanged. For $i = r-2$, we use the same rescalings, except that we
need to interpret the
second and third line as $$\widetilde \alpha_{r-2}=\alpha_{r-2}^{-1}.$$ In either case the rescaled tuple also belongs to $\mathcal C$.
Note also that although this rescaling may seem puzzling at first, its effect on the $x$ variables \eqref{xs} in the situation when $y_j
= \alpha_j,
\, 1\leq j \leq r-2$ is very simple: $x_j$ is unaffected for $j \neq i, r-1$, and the rescaled $x_i$ and $x_{r-1}$ are the old $x_{r-1}$ and $x_i$
respectively.
Recalling the definition of $\mathsf R$ in \eqref{R} we see that ${\mathsf R}$ is 
symmetric in $x_1, \ldots, x_{r-1}$, hence it is unmodified by the interchange of $x_i$ and $x_{r-1}$. Furthermore, under this interchange, 
$$\widetilde {y_1}^{m_1} \cdots \widetilde {y}_{r-1}^{m_{r-1}}=y_1^{m_1}\cdots y_{r-1}^{m_{r-1}}  \cdot (y_i \cdots y_{r-2})^{m_{i-1} -m_{i} + m_{r-1} - m_{r-2}}.$$
Since  
$$m_{i-1} -m_{i} + m_{r-1} - m_{r-2} \equiv 0 \mod N$$ by \eqref{reducedexp}, we deduce that
$$\Omega_1 (\alpha_1, \ldots, \alpha_{r-2}, y_{r-1}) = \Omega_1 (\widetilde \alpha_1, \ldots,\widetilde
\alpha_{r-2}, \widetilde y_{r-1}).$$
We therefore rewrite the sum over residues as $$- (r-1) \sum_{(\alpha_1, \ldots, \alpha_{r-2})\in \mathcal C}
{\text{Res}}_{y_{r-1} = 1} \Omega_1 (\alpha_1, \ldots, \alpha_{r-2}, y_{r-1}).$$
  
Repeating the procedure $r-1$ times we get
\begin{eqnarray}\label{itres} &&\sum_{(\alpha_1, \ldots, \alpha_{r-1}) \in {\mathcal C}} {\mathsf R} (\alpha_1 \cdot \ldots\cdot\alpha_{r-1} ,
\ldots, \alpha_{r-1} ,1) \cdot \alpha_1^{m_1}\cdots\alpha_{r-1}^{m_{r-1}}=\\ &=&(-1)^{r-1} (r-1)! \, \, {\res}_{y_1 = 1}
\ldots\res_{y_{r-1} = 1}  \Omega_{r-1} (y_1, \ldots, y_{r-1}),\nonumber \end{eqnarray} where $$\Omega_{r-1} = 
\frac{dy_1}{y_1}
\cdots \frac{dy_{r-1}}{y_{r-1}} \cdot\frac{N}{y_1^N -1} \cdots \frac{N}{y_{r-1}^N -1}\cdot \mathsf R(y_1 \cdots y_{r-1}, \ldots, y_{r-1} , 1) \cdot
y_1^{m_1}\cdots y_{r-1}^{m_{r-1}}.$$ 

In conclusion, collecting equations \eqref{vinow}-\eqref{itres}, we find that

\begin {equation}\label{intresa}\int_{\quot} \mathsf P(\bar \qa_2, \ldots, \bar \qa_r)\cdot 
\qa_r^{M}=\frac{(-1)^{\gbar\binom{r}{2}}}{r} N^{r\gbar+1} \cdot {\res}_{y_1 = 1} \ldots\res_{y_{r-1} = 1} \Omega_{r-1}.
\end {equation}

\section{Intersection numbers from $N$-asymptotics and Witten's sums} \label{der} \subsection {The intersections on the moduli
space of bundles.} To finish the proof of the main result, we use equation \eqref{finalint}. We need to evaluate the coefficient
of the highest power of $N$ in the expression \begin{equation}\label{intasres}\int_{\mathcal N}\exp({\bar f_2}) \mathsf P(\bar
a_2, \ldots, \bar a_r) =\left[N^{r(r-1)\gbar-\deg \mathsf P}\right] (-1)^{\gbar\binom{r}{2}}r^{\gbar} \,{\res}_{y_1 = 1}
\ldots\res_{y_{r-1} = 1} \Omega_{r-1}.\end{equation}

We first substitute $$y_i = \exp \left ( \frac{{\mathsf Y}_i}{N} \right )=\exp\left(\frac{\mathsf X_{i}-\mathsf
X_{i+1}}{N}\right), \, \, 1 \leq i \leq r-1,$$ where the variables $\mathsf X_1, \ldots, \mathsf X_r$ are defined by the system 
of
equations \eqref{xy} and \eqref{xy1}. Then, we rewrite the residue in \eqref{intasres} as \begin{eqnarray} \label{finalres} &\res_{\mathsf Y_{1} = 0}&\ldots \res_{\mathsf Y_{r-1} = 0} \frac{1}{e^{\mathsf Y_1} -1}\cdots\frac{1}{e^{ \mathsf Y_{r-1}} -1
}\cdot\mathsf R (e^{\frac{\mathsf X_1-\mathsf X_r}{N}}, \ldots, e^{\frac{\mathsf X_{r-1}-\mathsf X_r}{N}}, 1) \times \nonumber \\
&\times& \exp\left(\frac{m_1}{N}\mathsf Y_1+ \ldots + \frac{m_{r-1}}{N}\mathsf Y_{r-1} \right).\nonumber \end{eqnarray}

Note that since $\mathsf R$ is in fact a function of the normalized variables $\bar{x}_1, \ldots, \bar{x}_r,$ we have $$\mathsf R
(e^{\frac{\mathsf X_1-\mathsf X_r}{N}}, \ldots e^{\frac{\mathsf X_{r-1}-\mathsf X_r}{N}}, 1) = N^{r(r-1)\gbar-\deg \mathsf P}\, \mathsf
R( {\mathsf X_1}, \ldots, {\mathsf X_r}) + \, \text{lower order terms in} \, N.$$ Moreover, one sees from \eqref{min} that
$$\lim_{N\to\infty} \left(\frac{m_1}{N}\mathsf Y_1+\ldots+\frac{m_{r-1}}{N}\mathsf Y_{r-1}\right)=\left\{\frac{d}{r}\right\}\mathsf
Y_1+\ldots+\left\{\frac{d(r-1)}{r}\right\}\mathsf Y_{r-1}=\mathsf L.$$

Taking this into account, and using the definition of $\mathsf R$ in \eqref{Q} and \eqref{R}, we find that the highest order of $N$ has the coefficient
\begin{eqnarray}\label{finalfinalres}\res_{\mathsf Y_{1} = 0}\ldots \res_{\mathsf Y_{r-1} = 0} \frac{1}{e^{\mathsf Y_1}
-1} \cdots \frac{1}{e^{ \mathsf Y_{r-1}} -1 }\cdot \exp\left(\mathsf L \right)  \frac{ \mathsf P(s_2(\mathsf X_1, \ldots,
\mathsf X_r), \ldots, s_r(\mathsf X_1, \ldots, \mathsf X_r))}{\prod_{i<j}(\mathsf X_i-\mathsf X_j)^{2\gbar}}.\nonumber \end{eqnarray}

The statement of Theorem \ref{main} is now immediate in the light of equation \eqref{intasres}.

\begin {remark} One can include the odd cohomology classes $\bar b_{i}^{j}$ in the calculation by applying a suitable version of
Vafa-Intriligator involving odd $\mathsf b$ classes on $\quot$. Such statements are proved, in a particular case, in Proposition 
$2$ of \cite
{vi}, but the method used there extends in general. To keep the notation as simple as possible, we decided not to write down the general
formulas, leaving them to the interested reader. The remaining intersections involving other $f$ classes can be lifted easily to 
the Quot scheme and can in principle be computed there by equivariant localization as set up in \cite{vi}. 
We do not yet know of a systematic way of computing all intersections using the methods of this paper {\it i.e,} by studying 
$N$ asymptotics of the Vafa-Intriligator formula alone.   \end {remark}

\subsection{Witten's sums.}\label{irred} We rewrite the intersections computed in Theorem \ref{main} as infinite sums over the
representations of $SU(r)$, which is possible when the degree of $\mathsf P$ is small enough with respect to the genus $g$. We 
match the formulas
written down in \cite {witten1} \cite{liu}. Note that the arguments of Szenes-Jeffrey-Kirwan also show that the iterated residues
reproduce the Witten sums, but here we offer a direct derivation from equation \eqref{finalint}, bypassing the residue calculations of
Section \ref{residues}.

Using \eqref{finalint} and \eqref{vinow}, we equate $$\int_{\mathcal N}\exp({\bar f_2})\,\mathsf P(\bar a_2, \ldots, \bar a_r) $$ with the coefficient
of $N^{r(r-1)\gbar+1-\text {deg } \mathsf Q}$ in the sum over unordered tuples \begin{eqnarray}\frac{\mathsf u r^g}{r!} \sum_{\lambda_1, \ldots,
\lambda_r}
\frac{\left(\lambda_1\cdots\lambda_r\right)^{M-r\gbar}}{\prod_{i<j}\left(\left(\frac{\lambda_i}{\lambda_j}\right)^{\frac{1}{2}}-\left(\frac{\lambda_j}
{\lambda_i}\right)^{\frac{1}{2}}\right)^{2\gbar}}\cdot{\mathsf
Q}(\lambda_1,\ldots, \lambda_r).\nonumber\end{eqnarray} 
For each product $\lambda_1 \cdots \lambda_r$ there is a 
{\it unique} $\zeta$ such that $$\zeta^{r}=\lambda_1\cdots \lambda_r,$$ and such that, when setting $$\nu_i=\lambda_i\zeta^{-1}=\exp\left(\frac{2\pi
i}{N}\mu_i\right)$$ we have $$\mu_1+\ldots+\mu_r=0\,\,\, \text{with } -2N<\mu_i-\mu_j< 2N \text{ integers}.$$ Then, we compute
$$(\lambda_1\cdots\lambda_r)^{M-r\gbar}\cdot {\mathsf Q}(\lambda_1\cdots\lambda_r) =\zeta^{r(M-r\gbar)+\deg \mathsf Q} \cdot\mathsf Q(\nu_1, \ldots,
\nu_r)=\zeta^{N(d-r\gbar)} \cdot\mathsf Q(\nu_1, \ldots, \nu_r)$$ $$=\left(\frac{\lambda_r}{\nu_r}\right)^{N(d-r\gbar)}\cdot\mathsf Q(\nu_1, \ldots,
\nu_r)=\frac{1}{\nu_r^{Nd}}\cdot \mathsf Q(\nu_1, \ldots, \nu_r).$$ It remains to evaluate \begin{equation}\label{sinexp}\left[N^{r(r-1)\gbar-\deg
\mathsf Q}\right]\, \frac{(-1)^{d(r-1)}r^g}{2^{2\gbar\binom{r}{2}}r!}\,\sum_{\mu}\,\,\frac{\exp\left(-2\pi i d
\mu_r\right)}{\prod_{i<j}\left ( \sin \left( \frac{\mu_i-\mu_j}{N}\pi\right) \right )^{2\gbar}}\cdot {\mathsf Q}(e^{\frac{2\pi i}{N}\mu_1}, \ldots, 
e^{\frac{2\pi
i}{N}\mu_r}).\end{equation} A factor of $N$ disappeared in the normalization process to account for each possible value of the product $\lambda_1\cdots
\lambda_r$.

Now, \begin{equation}\label{error}\mathsf Q(e^{\frac{2\pi i}{N}\mu_1}, \ldots, e^{\frac{2\pi i}{N}\mu_r})=\frac{1}{N^{\deg \mathsf Q}}\cdot\mathsf Q(2\pi
i\mu_1, \ldots, 2\pi i\mu_r) +\text {lower order terms in } N\end{equation} where the coefficients of the lower order terms in $N$ are polynomials in $\mu$ of
smaller degree than $\deg \mathsf P$. Additionally, \begin{equation}\label{sin}\frac{1}{\sin^{2\gbar} x}=\frac{1}{x^{2\gbar}}+ \text{lower terms in }
x.\end{equation} Hence, we can evaluate the leading coefficient in \eqref{sinexp} to
\begin{equation}\label{sumsw}\frac{(-1)^{d(r-1)}r^g}{(2\pi)^{2\gbar\binom{r}{2}}r!}\cdot \sum_{\mu}\frac{\exp\left(-2\pi i d
\mu_r\right)}{\prod_{i<j}(\mu_i-\mu_j)^{2\gbar}} \cdot \mathsf Q(2\pi i\mu_1, \ldots, 2\pi i \mu_r).\end{equation}

The error terms in \eqref{error} and \eqref{sin} do not contribute. Indeed, the first order error terms uniformly dominate the rest, so it suffices to explain
that $$\frac{1}{N} \cdot \sum_{\mu} \frac{\text{polynomial in } \mu_1, \ldots, \mu_r}{\prod_{i<j} (\mu_i-\mu_j)^{2\gbar-2}}\to 0 \text{ as } N\to\infty.$$ This
is easy to see by expressing everything in terms of new variables $\sigma_i=\mu_i-\mu_{i+1}.$ We can uniquely solve for the $\mu$s using the constraint
$\mu_1+\ldots+\mu_r=0$. When $\deg \mathsf P$ is small compared to $g$, the degree of each $\sigma$ in the numerator is small compared to the degree of
$\sigma$ appearing in the denominator. In this case, the sum over $\mu$'s is convergent. The same argument gives the
convergence of the infinite sum \eqref{sinexp} when $\deg \mathsf P$ is small.

We now express the result in terms of the representation theory of $SU(r)$. We write $e_i$ for the coordinates on the dual Cartan algebra, and agree that
$e_i-e_j$, $i<j$ are the positive roots of $SU(r)$. Let $$\rho=\sum_{i} \frac{r-2i+1}{2} e_i$$ be half the sum of the positive roots. We order $\mu_1>
\mu_2>\ldots>\mu_r.$ Setting $$\chi_i=\mu_i-\frac{r-2i+1}{2},$$ we have $$\chi_1\geq \ldots\geq \chi_r, \,\,\chi_i-\chi_j\in \mathbb Z, \,\,\,
\chi_1+\ldots+\chi_r=0.$$ Thus, we can think of $$\chi=\chi_1 e_1+\ldots+\chi_r e_r$$ as the highest weight of an irreducible representation $\mathcal
R_{\chi}$ of $SU(r)$. The Weyl dimension formula gives $$\text{dim }\mathcal
R_{\chi}=\prod_{i<j}\frac{\chi_i-\chi_j+j-i}{i-j}=\prod_{i<j}\frac{\mu_i-\mu_j}{i-j}.$$ Moreover, the scalar action of the central element
$$c=\exp\left(\frac{2\pi i d}{r}\right) I$$ has $\chi$-trace $${\text {Trace}_{\chi}(c)}=\exp(-2\pi i d \chi_r)\cdot {\text{dim }\mathcal
R_{\chi}}=(-1)^{d(r-1)}\exp(-2\pi i d\mu_r)\cdot {\text{dim }\mathcal R_{\chi}}.$$

Putting everything together, we see that \eqref{sumsw} transforms into the Witten sums \begin{equation}\label{volume}\int_{\mathcal N}\exp({\bar f_2})\mathsf
P(\bar a_2, \ldots, \bar a_r)=\mathsf C\cdot\sum_{\chi}\frac{\text {Trace}_{\chi}(c)}{\left(\text {dim } \mathcal R_{\chi}\right)^{2g-1}}\cdot
\mathsf Q(2\pi i (\chi+\rho))\end{equation} when $\deg \mathsf P$ is small compared to the genus. Here, $\mathsf Q$ is defined in \eqref{Q}, and we used the constant $$\mathsf
C=\frac{r^{g}}{(2\pi)^{r(r-1)\gbar}\cdot 1!^{2\gbar}\cdots(r-1)!^{2\gbar}}.$$

\section{Vanishing of intersections on the Quot scheme.} 

In this final section, we aim to establish the following vanishing statement

\begin{proposition} \label{vanishing} Assume that $r$ and $d$ are relatively prime, and $r\geq 2$. Let
$\mathsf P(\bar \qa_2, \ldots,
\bar \qa_r)$ and $\mathsf S(\qa_1, \qa_2, \ldots, \qa_r)$ be polynomials on $\quot$, such that the weighted
degree $\deg \mathsf P > r(r-1)\gbar$, so that $\deg {\mathsf P} + \deg {\mathsf S} < \frac{N}{r}$, and so that
$\deg {\mathsf P} + \deg {\mathsf S}
+ rM = Nd-r(N-r) \gbar$ for a positive integer $M$. Then we have \begin{equation} \label{weak0} \int_{\left[\quot\right]^{vir}}
\mathsf P(\bar \qa_2, \ldots,
\bar\qa_r)\cdot \mathsf S(\qa_1, \qa_2, \ldots, \qa_r) \cdot \qa_r^M = 0.\end{equation}

\end{proposition}

Note that Theorem $5$ of \cite{vi} implies that the bound on the degree of $\mathsf P$ {\it cannot} be lowered. 

Proposition \ref{vanishing} gives constraints governing the virtual number of maps from $C$ to the Grassmannian with incidence
conditions to special Schubert subvarieties at fixed domain points; as mentioned in the introduction, these numbers are actual counts provided that
the degree $d$ is large \cite {bertram}. The exact equations are obtained by linearity from the $a_i$-monomials of the product
$\mathsf P\cdot \mathsf S$, by requiring incidences at distinct domain points for each occurrence of $a_i$. It would be interesting to interpret these
enumerative constraints geometrically. \vskip.1in
{\bf Proof.} Fix polynomials $\mathsf P$ and $\mathsf S$ as in the statement of the proposition. Using the Vafa-Intriligator formula we obtain the following sum over unordered tuples of distinct roots
of unity $$\int_{\left[\quot\right]^{vir}}\mathsf P(\bar\qa_2, \ldots, \bar \qa_r)\cdot \mathsf S(\qa_1, \ldots,
\qa_r) \cdot \qa_r^M = $$ $$ = \frac{\mathsf u N^{r\gbar}}{r!} \sum_{\lambda_1, \ldots, \lambda_r} \mathsf R(\lambda_1, 
\ldots, \lambda_r)\cdot {{\mathsf 
T}}(\lambda_1,
\ldots, \lambda_r) \cdot \left(\lambda_1\cdots\lambda_r\right)^{M-\gbar},$$ where ${{\mathsf T}}$ is the polynomial $\mathsf S$ expressed in terms of the Chern roots, 
and $\mathsf R$ is
defined as before by equations \eqref{Q} and \eqref{R}. We rescale the variables to obtain $\lambda_r=1$, and
express everything in terms of the $\alpha$s defined in \eqref{alpha}. We are then to prove that 
\begin{equation}\label{weakvan}\nonumber \sum_{(\alpha_1, \ldots,
\alpha_{r-1})\in \mathcal C}\mathsf R(\alpha_1\cdot\ldots\cdot\alpha_{r-1},\ldots, \alpha_{r-1}, 1) \cdot
{{\mathsf T}} (\alpha_1\cdot\ldots\cdot\alpha_{r-1},\ldots, \alpha_{r-1}, 1) \cdot \alpha_1^{m_1}\cdots\alpha_{r-1}^{m_{r-1}}=0.\end{equation} 
Here $m_1, \ldots, m_{r-1}$ are defined as in \eqref{reducedexp}. By the assumption on the degrees of $\mathsf R$ 
and 
$\mathsf S$, the meromorphic form 
\begin{eqnarray} \Omega_{r-1}& = &\frac{dy_1}{y_1}
\cdots \frac{dy_{r-1}}{y_{r-1}} \cdot\frac{N}{y_1^N -1} \cdots \frac{N}{y_{r-1}^N -1}\cdot \mathsf R(y_1 \cdots y_{r-1}, \ldots, y_{r-1} , 1) \cdot
\nonumber \\
 & \cdot & {{\mathsf T}} (y_1 \cdots y_{r-1}, \ldots, y_{r-1} , 1) \cdot 
 y_1^{m_1}\cdots y_{r-1}^{m_{r-1}} \nonumber
\end{eqnarray}
associated with the above sum has 
poles only at the 
$N^{\text{th}}$ roots of $1$. Thus, by \eqref{itres}, it suffices to show
that $$\res_{y_1 = 1} \ldots \res_{y_{r-1} = 1} \Omega_{r-1}(y_1, \ldots, y_{r-1})=0.$$ To prove this, we make use of a trick that we learned from \cite {EK}. Using the auxiliary variable $t$, we substitute
$$y_i=\exp\left(\frac{t\mathsf Y_i}{N}\right)=\exp\left(\frac{t\left(\mathsf X_i-\mathsf X_{i+1}\right)}{N}\right).$$ We need to show the vanishing
of the iterated residue \begin{eqnarray} \label{tres} &\res_{\mathsf Y_{1} = 0}& \ldots \res_{\mathsf Y_{r-1} = 0}
\frac{t}{e^{t\mathsf Y_1} -1}\cdots\frac{t}{e^{t \mathsf Y_{r-1}} -1 }\cdot\mathsf R \left(e^{\frac{t(\mathsf X_1-\mathsf X_r)}{N}}, \ldots,
e^{\frac{t(\mathsf X_{r-1}-\mathsf X_r)}{N}}, 1\right) \cdot \nonumber \\ & \cdot& {{\mathsf T}} 
\left(e^{\frac{t(\mathsf X_1-\mathsf X_r)}{N}}, \ldots,
e^{\frac{t(\mathsf X_{r-1}-\mathsf X_r)}{N}}, 1\right) \cdot \exp\left(\frac{m_1}{N}\cdot t\mathsf Y_1+ \ldots +
\frac{m_{r-1}}{N} \cdot t \mathsf Y_{r-1} \right) .\end{eqnarray} We have already observed in Section \ref{der} that $$\mathsf R 
\left(e^{\frac{t(\mathsf X_1-\mathsf X_r)}{N}},
\ldots, e^{\frac{t(\mathsf X_{r-1}-\mathsf X_r)}{N}}, 1\right)=\left(\frac{t}{N}\right)^{\text {deg }\mathsf P-r(r-1)\gbar} \mathsf R(\mathsf X_1,
\ldots, \mathsf X_r)+\text {higher terms in } t.$$ All other terms in \eqref{tres} are holomorphic in $t$, hence the order in $t$ of the residue
above is at least $\deg \mathsf P-r(r-1)\gbar\geq 1$. On the other hand, the intersection number this residue computes is independent of $t$, and
therefore it must vanish.\qed

\begin{remark} By transferring the intersection \eqref{weak0} to $\p_{N, r,d}$ and pushing it forward to 
$\mathcal M$, Proposition \ref{vanishing} is seen to follow from the vanishing of the Pontryagin ring of 
$\mathcal M$ in degree greater than $r(r-1)\gbar$, proved in \cite {EK}. 

Conversely, it is tempting to speculate that 
the vanishing statement \eqref{weak0} implies the vanishing of the Pontryagin ring of $\mathcal N.$ This is true when $r=3$ as shown below.
In general, taking $\mathsf S=1$ and considering the highest order term in $N$ of the pushforward along $\pi$ yields $$\int_{\mathcal N}\exp (\bar f_2) \cdot \mathsf P(\bar a_2, \ldots, \bar a_r)=0 \,\,\text{ for }\deg\mathsf P>r(r-1)\gbar.$$ An analysis of the full $N$ asymptotics of the pushforward of $\mathsf
S(\qa_1, \ldots, \qa_r) \cdot \qa_r^{M}$ may be possible, if cumbersome. To start, one needs to examine 
the full $N$-asymptotics of the Segre class of Lemma \ref{segreasymp}, given by 
$$s\left(\mathcal H^N\right)=\exp\left(N\cdot \sum_{j}(-1)^{j+1}j!\,\text{ch}_{j+1}(\mathcal
H ) \right).$$ Now take $r=3$ and $\mathsf S=1$. Using the above expression for the Segre class, and observing that $\text {ch}_{j+1}(\mathcal H)$ is linear in $\bar f_3$ when $j\geq 1$, we see that the pushforward $\pi_{\star} \left(a_3^M\right) \mathsf P(\bar a_2,
\bar a_3)$ on $\mathcal N\times J$ equals $$\sum_{k} \frac{1}{r^{2g}}\cdot \frac{N^{e(\mathsf P)-k}}{k!} \cdot \left(\bar f_3^{k}\, \exp(\bar f_2+ r\Theta)\, \mathsf P(\bar a_2, \bar a_3) + \text {lower terms in }
\bar f_3\right).$$ Using induction on $k$, we obtain the
vanishing of all intersection products $$\int_{\mathcal N} \bar f_3^{k}\, \exp({\bar f_2})\, \mathsf P(\bar a_2,
\bar a_3)=0.$$ The rank $3$ Pontryagin vanishing follows, upon including the odd $b$ classes in the Vafa-Intriligator
formula. In higher rank, it is to be surmised that the consideration of the full $N$ asymptotics of the pushforwards $\mathsf
S(\qa_1, \ldots, \qa_r) \cdot \qa_r^{M}$ for {\em{all}} polynomials $\mathsf S$ subject to the assumption of 
Proposition \ref{vanishing}, will inductively give the vanishing of all evaluations 
$$\int_{\mathcal N} {\bar{f}}_2^{k_2} \cdots {\bar{f}}_{r}^{k_r} \cdot \mathsf P(\bar{a}_2, \ldots, \bar{a}_{r} ).$$ \end{remark}

\end{document}